\documentclass[12pt]{article}
\usepackage{amsmath,amsfonts,amssymb}
\def\N{\mathbb N}

\def\Z{\mathbb Z}
\def\C{\mathbb C}
\def\Q{\mathbb Q}

\def\q{\mathfrak q}

\def\G{\mathbb G}
\def\ZZ{\mathfrak Z}
\def\SLL{{\rm SLL}}
\def\GLL{{\rm GLL}}

\def\GA{{\rm I}\!\Gamma}
\def\fin{\hbox{\quad $\blacksquare$}}
\def\Cal#1{{\cal #1}}

\def\demo{{\it Proof.\quad  }}
\def\pmatrix#1{\left[\begin{matrix}#1\end{matrix}\right]}

\author{Carmen Rosa Giraldo Vergara \\
\footnotesize Instituto de Matem\'atica, Universidade Federal do Rio de Janeiro\\
\footnotesize CEP 21945-970, Rio de Janeiro, Brasil\\
\footnotesize E-mail: carmita@pg.im.ufrj.br\\ \and
Fabio Enrique Brochero Mart\'{\i}nez\\
\footnotesize Instituto de Matem\'atica Pura e Aplicada\\
\footnotesize  CEP 22460-320,  Rio de Janeiro RJ, Brasil \\
\footnotesize E-mail: fbrocher@impa.br }
\date{
}
\title{Zorn's Matrices and finite index subloops\footnote{ This research was supported by CNPq, Brasil. }}
\newtheorem{theorem}{Theorem}[section]
\newtheorem{definition}{Definition}[section]
\newtheorem{lema}{Lemma}[section]
\newtheorem{corolario}{Corollary}[section]
\newtheorem{proposicion}{Proposition}[section]
\newtheorem{Observacion}{Remark}[section]

\begin{document}

 \maketitle
\begin{abstract}
The Zorn's Algebra  $\ZZ(R) $ has a multiplicative function called
determinant with properties similar to the usual one. The set of
elements in $\ZZ(R)$ with determinant 1 is a Moufang loop that we
will denote by $\GA$. In our main result we prove that if $R$ is a
Dedekind algebraic number domain that contains an infinite order
unit, each finite index  subloop $\Cal L$, such that $\GA$ has the
weak Lagrange property relative to $\Cal L$, is congruence
subloop. In addition, if $R=\Z$, then we present normal subloops
of finite index in $\GA$ that are not congruence subloops.
\end{abstract}
\section{Introduction}
Let $R$ be a commutative ring with unit 1, $I$ an ideal of  $R$
and $SL(n,R)$
 the $n{\times}n$ special linear group over $R.$
The principal congruence group of level
 $I$ in $SL(n,R)$ is a set of  matrices congruent to the identity modulo the ideal $I.$
 It said that
$SL(n,R)$ satisfies the congruence subgroups property if every
finite index subgroup contains a principal congruence group.

In the last decades the congruence groups  achieved  own
relevance, different from the traditional application in the
geometric field about the classification of elliptic curves over
$\C $ and the study of modular forms. This relevance is due to the
works of Mennicke, Serre, Lazard, Bass, Vaser\v{s}te\v{\i}n,
Newman and some others. Mennicke \cite{Men} proved  that
$SL(n,\Z)$ with $n\geq 3$ satisfies the congruence subgroups
property. In addition Bass, Milnor, Serre \cite{BMS} proved that
$SL(n,R)$ satisfies that property for $n\geq3$ for an ample
variety of rings (in particular for any ring o algebraic numbers).
For $n=2$, Wohlfahrt (See \cite{Rei}) showed a criterion of
determing when a finite index subgroup is a congruence group, and
used this criterion to show that in $SL(2,\Z)$ exists finite index
subgroups that is not congruence subgroup. In general,
\begin{theorem}[Serre \cite{Ser}] Let $R$ be an  algebraic
integer domain that contains an infinite order unit. Then
$SL(2,R)$ satisfies the congruence subgroups property.
\end{theorem}
 In that article, Serre also proved that if $R$ is
an algebraic integer ring $\Cal O_d$ of the field $\Q(\sqrt{-d})$
with $d\in \N^*$, then the group $SL(2,\Cal O_d)$ does not have
the congruence subgroups property.

In addition, Varser\v{s}te\v{\i}n  showed a relationship between
congruence subgroup and groups generates by elementary matrices
when $R$ is a Dedekind algebraic domain that contains infinitely
many units.

\begin{theorem}[Varser\v{s}te\v{\i}n \cite{Var}] Let $R$ be a Dedekind
algebraic numbers  domain that contains  infinitely many units,
$\q$ an ideal of $R$ and $E(\q)$ the subgruop generates by the
matrices $\left(\begin{array}{cc} 1&a\\0&1\end{array}\right)$,
$\left(\begin{array}{cc} 1&0\\a&1\end{array}\right)$ with $a\in
\q$. Then the group $E(\q)$ has finite index in $SL(2,R)$, in
particular $E(R)=SL(2,R)$.
\end{theorem}

On the other side,  Zorn's algebra $\ZZ(R)$  contains Moufang
loops analogous to the groups $GL(2,R)$ and $SL(2,R)$. In fact,
denote by $R^3$ the three dimensional  vector space over $R$.
Zorn's Algebra $\ZZ(R)$ over  $R$ is the set of $2{\times}2$ matrices
$$\left[\begin{array}{cc} a& x\\ y & b\end{array}\right]\quad
a,b\in R \quad x,y\in R^3,$$ with the binary operations sum and
product, where the sum is defined by the natural form, component
to component, and the product is given by the rule
$$\left[\begin{array}{cc} a_1& x_1\\ y_1 &
b_1\end{array}\right]\left[\begin{array}{cc} a_2& x_2\\ y_2 &
b_2\end{array}\right]=\left[\begin{array}{cc} a_1a_2+x_1{\cdot}y_2&
a_1x_2+b_2x_1-y_1{\times}y_2\\ a_2y_1+b_1y_2+x_1{\times}x_2 &
b_1b_2+y_1{\cdot}x_2\end{array}\right]$$ where  $\cdot$ and $\times$
denote the dot and cross vectorial product in $R^3$.

\noindent The determinant function $\det:\ZZ(R)\to R$  defined by
\hbox{$\det(A) =ab-x{\cdot}y$}, where $A=\left[\begin{array}{cc} a& x\\
y & b\end{array}\right]$, is a multiplicative function . Thus, an
element $A$ is invertible if and only if $\det(A)\in R^*$ and
then \hbox{$A^{-1}=\frac 1{\det A}\left[\begin{array}{cc}b& -x\\
-y & a\end{array}\right]$}. Zorn's Algebra is alternative, it
follows that the invertible elements set is a Moufang Loop. This
set is called a general linear loop and it is denoted by
$$ \GLL(2,R)=\{A\in \ZZ(R)|\ \det (A)\in R^*\}.$$
Similarly, we define the special linear loop as follows
$$ \GA= \SLL(2,R)=\{A\in \GLL(2,R)| \ \det A= 1\}.$$

We are going to developed for  loops an analogue theory  to the
congruences groups theory, where $R$ is a Dedekind algebraic
numbers domain. In particular, if $R$  contains an infinite order
unit, we will  prove an analogous Serre's Theorem for these loop.
In addition, if $R=\Z$, we find a family of finite index subloops
that are not  congruence loops.

\section{Congruence subloop and finite index subloop}

 Let $\q$ be an ideal of $R$. We define a principal congruence
subloop of $\GA$ of level  $\q$ as a set of all matrices $A$ of
$\GA$ such that
$$A\equiv I \pmod \q,$$ where the congruence is component by component. This loop is denoted by  $\GA(\q)$.
In particular, $\GA(R)=\GA$. A subloop of  $\GA$ is called a
congruence subloop if it contains a principal congruence subloop
$\GA(\q)$ for some ideal $\q$ of $R$.

In the follow, $R$ denote a Dedekind  domain, and for each $\q$
non zero  ideal of $R$,  $\Delta(\q)$ denotes the smallest normal
subloop of
 $\GA$ that contains every matrix of the form
$\left[\begin
{array}{cc}1 & x \\
0 & 1\end{array}\right]$, $\left[\begin
{array}{cc}1 & 0 \\
y & 1\end{array}\right]$ where $x,y \in \q^3.$

\begin{definition} Let  $\Cal L$ be a  subloop of $\GA$, and suppose that the set
 $$S=\{\q\subset R | \Delta (\q)\subset \Cal L\}$$ is not empty.
Define the level of   $\Cal L$ as the maximal element of $S$.
\end{definition}

\begin{lema}\label{fundamental} Let  $\q$ be an non zero ideal of $R$ and $A\in \GA(\q)$. If $A=\left[\begin
{array}{cc}1 & v \\
u & b\end{array}\right]$, then $A$ can be written as product of
the matrices $\left[\begin
{array}{cc}1 & ae_j \\
0 & 1\end{array}\right]$ and  $\left[\begin
{array}{cc}1 & 0 \\
ae_j & 1\end{array}\right]$, where $a\in \q$,  $e_1=(1,0,0)$,
$e_2=(0,1,0)$ and $e_3=(0,0,1)$.
\end{lema}
\demo
First, suppose $A$ has the form   {\small $\left[\begin{array}{cc} 1&(0,0,0)\\
(u_1,u_2,u_3)& 1\end{array}\right]$} then define
 $B=\left[\begin{array}{cc} 1&(0,0,0))\\
(0,u_2,u_3)& 1\end{array}\right]\in \GA(\q)$, thus
$$\begin{array}{rcl} A&=&B+\left[\begin{array}{cc} 0&0\\
u_1e_1& 0\end{array}\right]\\ &=&
B\left(I+B^{-1}\left[\begin{array}{cc} 0&0\\ u_1e_1&
0\end{array}\right]\right)=B\left[\begin{array}{cc}
1&(0,u_1u_3,-u_1u_2)\\ (u_1,0,0)& 1\end{array}\right]\\ &=&
B\left(\left[\begin{array}{cc} 1&0\\(u_1,0,0) &
1\end{array}\right]\left[\begin{array}{cc} 1&(0,u_1u_3,-u_1u_2)\\
0& 1\end{array}\right]\right)
\end{array}
$$
By the same procedure, we obtain that $$B=\left[\begin{array}{cc}
1&0\\u_3e_3 & 1\end{array}\right]\left(\left[\begin{array}{cc}
1&-u_2u_3e_1\\0 & 1\end{array}\right]\left[\begin{array}{cc}
1&0\\u_2e_2 & 1\end{array}\right]\right)$$ and
$$\left[\begin{array}{cc} 1&(0,u_1u_3,-u_1u_2)\\
0& 1\end{array}\right]=\left[\begin{array}{cc} 1&-u_1u_2e_3\\0 &
1\end{array}\right]\left(\left[\begin{array}{cc}
1&0\\u_1^2u_2u_3e_1 & 1\end{array}\right]\left[\begin{array}{cc}
1&u_1u_3e_2\\0 & 1\end{array}\right]\right)$$

\noindent In the general case, if   $A=\left[\begin
{array}{cc}1 & (v_1,v_2,v_3) \\
(u_1,u_2,u_3) & b\end{array}\right]$, we have
$A=(((CA_3)A_2)A_1$ where  $A_j=\left [\begin {array}{cc} 1&v_je_j\\
0&1\end {array}\right]$ and
$ C= \left [\begin {array}{cc} 1&(0,0,0)\\
(u_1+v_3v_2,u_2-v_3v_1,u_3+v_2 v_1)&1\end {array}\right ]. $  The
result follows from the first case.\fin

Let  $\GA_{(j)}(\q)$ denote a  subloop {\small
$$\GA_{(j)}(\q):=\left\{ \pmatrix{a&be_j\cr ce_j&d}\in \GA\left|
\pmatrix{a&b\cr c&d}\equiv I
\pmod \q\right.\right\}\quad\hbox{and}\quad
GA_{(j)}:=GA_{(j)}(R).$$} We are going to show that these three
subloops generate the loop $\GA(\q)$, but before we need a result
from commutative ring theory.
 Let  $R$ be  a commutative
ring with unit  1, and $M$ be a $R$-module. $R$ is called a local
ring if there is an unique local maximal ideal, and when there are
a finite number of maximal ideals, it is called semilocal.  An
element  $x\in M$ is called  unimodular  if there is a linear form
$L:M\to R $ such that $L(x)=1$.

\begin{lema}\label{unimodular} Let  $x=(x_1,\dots, x_m)$ be an unimodular element of
 $R^m$. If  $R$ is a semilocal noetherian ring, then there are  $y_2,\dots, y_m\in R$
 such that
$$x_1+y_2x_2+\cdots+y_mx_m$$ is invertible in  $R$.
\end{lema}
\demo  See \cite{BLS} page 386.\fin

\begin{theorem}\label{genera-n} Let $R$ be a Dedekind domain. For every $\q$ non zero  ideal of $R$,
$\GA(\q)$ is generated by $\GA_{(1)}(\q)$, $\GA_{(2)}(\q)$ and
$\GA_{(3)}(\q)$.
\end{theorem}
\demo Let  $\Cal L=\langle
\GA_{(1)}(\q),\GA_{(2)}(\q),\GA_{(3)}(\q)\rangle$. It is clear
that
 $\Cal L\subset \GA(\q)$, thus, we only need to show  $\GA(\q)\subset \Cal L$. Let  $A\in \GA(\q),$
 i.e.,
   $A$  is  a matrix of the form
$\pmatrix{a&(u_1,u_2,u_3) \cr (v_1,v_2,v_3)&b}$ with
$ab-(v_1u_1+v_2u_2+v_3u_3)=1$ and $a$ and $b$   congruent to 1
modulo $q$, in particular,   $a\ne 0$. If $\q=R$ and
 $a=0$ , since  $ab-(v_1u_1+v_2u_2+v_3u_3)=1$, there is  $j$
such that $v_ju_j\ne 0$.
It follows that it is sufficient to prove   $$T_jA=\left[\begin{array}{cc} e_j{\cdot}u&be_j+e_j{\times}u\\
-ae_j+e_j{\times}v& -e_j{\cdot}v\end{array}\right]$$ is in  $\Cal L$, and
therefore we can suppose, in this case too, that $a\ne 0$.

Since $R$ is Dedekind, $\frac R{(a)}$ is a semilocal ring and
$(-v_1)u_1+(-u_2)v_2+(-u_3)v_3\equiv 1 \pmod a$,
 thus $(u_1,v_2,v_3)$ is  unimodular in  $\left(\frac R{(a)}\right)^3$. Therefore from the  lemma
\ref{unimodular}, exist $t$ and $s$ such that $u_1+v_2t+v_3s=u_1'$
is invertible in  $\frac R{(a)}$ and it follows that $a$ and
$u'_1$ are relative primes. Define
$$u_2'=-v_1t+u_2\quad\hbox{and}\quad u_3'=-v_1s+u_3.$$
It is easy to prove that $B=\pmatrix{a&(u'_1,u'_2,u'_3) \cr
(v_1,v_2,v_3)&b}\in \GA(\q)$ and $$A=B\pmatrix{1&
b(u_1-u'_1,u_2-u'_2,u_3-u'_3) \cr
-(u'_2u_3-u'_3u_2,u'_3u_1-u'_1u_3,u'_1u_2-u'_2u_1)&1}.$$ It
follows from  lemma \ref{fundamental} that we only need to prove
that $B\in \Cal L$. Let $q$ be an arbitrary element of $\q^*$.
Since \hbox{$(a,qu'_1)=1$}, there are $x,y$ integer numbers such
that $ax+qu'_1y=1$. But $a\equiv 1\pmod \q$, therefore $x\equiv
1\pmod \q$ and $\left[\begin{array}{cc}
a&(u'_1,0,0)\\(-qy,0,0)&x\end{array}\right]\in \GA_{(1)}(\q)$.
Thus

{\scriptsize$$\displaylines{\left[\begin{array}{cc} a&(u'_1,u'_2,u'_3)\\
(v_1,v_2,v_3)& b\end{array}\right]=\hfill\cr =\left [\begin
{array}{cc}
1&(0,au'_2-v_3qy,au'_3+v_2qy)\\
(xv_1+bqy,xv_2-u'_3u'_1,xv_3+u'_2u'_1)&ba- v_1u'_1\end
{array}\right ]\left[\begin{array}{cc}
a&(u'_1,0,0)\\(-qy,0,0)&x\end{array}\right]}$$} Finally,  lemma
\ref{fundamental} shows that  $B\in \Cal L$.\fin

\begin{corolario} Let $R$ be a Dedekind algebraic domain that
contains  infinitely many units. Then $\GA$ is generated by the
matrices $\left[\begin
{array}{cc}1 & ae_j \\
0 & 1\end{array}\right]$ and  $\left[\begin
{array}{cc}1 & 0 \\
ae_j & 1\end{array}\right]$ where $a\in R$ and $j=1,2,3$.
\end{corolario}
\demo It follows from Vaser\v{s}te\v{\i}n theorem.\fin

\noindent The following theorem  is a generalization of
Wohlfahrt's criterion to $\GA.$
\begin{theorem}
Let  $\q_1,\q_2$ be ideals of the Dedekind domain $R$.  Then
$\GA(\q_1)\subseteq \Delta(\q_1)\GA(\q_2)$.
\end{theorem}
\demo This proof is similar to the proof of   Wohlfahrt's theorem
 made by Mason and Stothers (see \cite{MsSt}).
Let  $A=\left[\begin
{array}{cc}a & (v_1,v_2,v_3) \\
(u_1,u_2,u_3) & b\end{array}\right]$ be an arbitrary element of
$\GA(\q_1)$, we need to find $B\in \Delta(\q_1)$ such that
$A\equiv B\pmod {\q_2}$, and therefore  $A=B(B^{-1}A)$ where
$B^{-1}A\in \GA(\q_2)$.

\noindent{\it Case 1:} If $a\equiv 1 \pmod {\q_2}$, then it is
sufficient to take {\small $B=\left[\begin
{array}{cc}1 & (v_1,v_2,v_3) \\
(u_1,u_2,u_3) & ab\end{array}\right]$}, because    lemma
\ref{fundamental} shows that   $B\in \Delta(\q_1).$ In addition,
\begin{eqnarray*}B^{-1}A&=&\left[\begin
{array}{cc}a^2b-(u_1v_1+u_2v_2+u_3v_3) & (ab-b)(v_1,v_2,v_3) \\
(1-a)(u_1,u_2,u_3) & b-(u_1v_1+u_2v_2+u_3v_3)
\end{array}\right]
\\
&=&I+(a-1)\left[\begin
{array}{cc}ab& b(v_1,v_2,v_3) \\
-(u_1,u_2,u_3) & -b\end{array}\right] \equiv I\pmod
{\q_2}.\end{eqnarray*}
 Thus  $A\in \Delta(\q_1)\GA(\q_2).$

\noindent{\it Case 2:} If  $(a,\q_2)=1$, then the congruence
$ax\equiv 1 \pmod {\q_2}$ has solution. Let  $a'$ be a solution,
$c=gcd(v_1,v_2,v_3)$, $v'_j=\frac {v_j}c$ for $j=1,2,3$ and
$X=\left[\begin
{array}{cc}1 & a'(1-a-c)(v'_1,v'_2,v'_3) \\
(0,0,0) & 1\end{array}\right]$. Notice that $a'(1-a-c)$ is in
$\q_1$ since $1-a\in \q_1$  and $c\in \q_1$,  thus  $X\in
\Delta(\q_1).$ In addition,
$$AX\equiv \left[\begin
{array}{cc}a & (1-a)(v'_1,v'_2,v'_3) \\
(u_1,u_2,u_3) & b+a'( 1-a-c)(\frac{ab-1}c)\end{array}\right]\pmod
{\q_2}.$$ Define $T_1=\left[\begin
{array}{cc}1 & (0,0,0) \\
(t_1,t_2,t_3) & 1\end{array}\right]$ where
$v'_1t_1+v'_2t_2+v'_2t_2=1$. We have
$$T_1^{-1}(AX)T_1\equiv \left[\begin
{array}{cc}1 & *
\\
\ *  & *\end{array}\right]\pmod {\q_2}.$$ It follows from  {\it
Case 1} that $T_1^{-1}(AX)T_1\equiv B\pmod {\q_2}$ for some $B\in
\Delta(\q_1)$, therefore $A\equiv (T_1BT_1^{-1})X^{-1}\pmod
{\q_2}$, with $(T_1BT_1^{-1})X^{-1}\in \Delta(\q_1)$.

\noindent{\it Case 3:} In the general case,  denoting
$d=(u_1v_1+u_2v_2+u_3v_3)$, we have  $ab-d=1$.  Then $(a,d)$ seen
as an element of   $\left(\frac R{\q_2}\right)^2$ is unimodular,
and  since  $\frac R{\q_2}$ is semilocal ring, from lemma
\ref{unimodular}, exists $t$ such that $a-td$ is invertible in
$\frac R{\q_2}$. Define $T=\left[\begin
{array}{cc}1 & (0,0,0) \\
-t(u_1,u_2,u_3) & 1\end{array}\right]$, then
$$T^{-1}AT=\left[\begin
{array}{cc}a-td & (v_1,v_2,v_3) \\
-(-at-1+bt+{t}^{2}d )(u_1,u_2,u_3) & b+td\end{array}\right],$$
where  $(a+td,\q_2)=1$, and from  {\it Case 2}, $T^{-1}AT\equiv
B\pmod {\q_2}$ for some $B\in \Delta(\q_1)$. Therefore $A\equiv
TBT^{-1}\pmod {\q_2}$, with $TBT^{-1}\in \Delta(\q_1)$.\fin

\begin{corolario}\label{Wolf} Let  $\Cal L$ be a congruence subloop of $\GA$ of level  $\q$.
Then $\Cal L\supset \GA(\q)$.
\end{corolario}
\demo Since  $\Cal L$ is a congruence subloop, then there is an
ideal  $\q'$ such that  $\GA(\q')\subset\Cal L$, furthermore
$\Delta(\q)\subset \Cal L$, then  $\GA(\q)\subset
\Delta(\q)\GA(\q')\subset \Cal L$.\fin

To show  a generalization of Serre's Theorem to $\GA$, we need the
following fact about Loop Theory

\begin{definition} Let $\Cal L$ be  a loop with the inverse property and
$\Cal H$ a subloop of $\Cal L$. We said that $\cal L$ has the
Lagrange property relative to $\Cal H$,  if $\Cal L$ and $\Cal H$
satisfy one of the following equivalent conditions (see \cite{GJM}
page 52)
\begin{enumerate}
\item $\Cal Hx\cap \Cal Hy\ne \emptyset$, $x,y\in \Cal L$ if and
only if $\Cal Hx=\Cal Hy$.
\item $\Cal H(hx)=\Cal Hx$ for all $x\in \Cal L$ and $h\in \Cal H$
\end{enumerate}
We said that $\cal L$ has the weak Lagrange property relative to
$\Cal H$, if there exists $\Cal F$, finite index subloop of $\Cal
H$, such that $\Cal L$ has the Lagrange property relative to $\Cal
F$.
\end{definition}
Observe that if $\Cal H$ is normal subloop of $\Cal L$, then $\Cal
L$  has the  Lagrange property relative to $\Cal H$, and there
exist subloops, such that $\GA$  have not the Lagrange property
relative to these subloops. For instance,
$$\GA_{(1)}\pmatrix{0&e_1+e_2\cr -e_2+e_3&0}\subset\left\{\pmatrix
{0&(a,b,c)\cr (0,d,e)&c}\in \GA\right\}$$ and {\footnotesize
$$\pmatrix{1&e_1\cr e_1&2}\left(\pmatrix {1&0\cr
e_1&1}\pmatrix{0&e_1+e_2\cr
-e_2+e_3&0}\right)=\pmatrix{0&(2,3,2)\cr (0,-3,4)&3}\notin
\GA_{(1)}\pmatrix{1&e_1+e_2\cr -e_2+e_3&0},$$}%
thus $\GA$ has not the Lagrange property  relative to $GA_{(1)}$.
Similarly, $\GA$ has not the Lagrange property  relative to
$\GA_{(1)} \GA(n)$, for all $n\ge 2$.
 In addition, when
$\Cal L$ has the Lagrange property relative to $\Cal H$, there
exists a subset $T$ of $\Cal L$ called a transversal of $\Cal H$
such that
$$\Cal L=\bigcup_{x\in T} \Cal Hx\quad\hbox{with} \quad \Cal Hx\cap
\Cal H y=\emptyset\hbox{ for } x,y\in T.$$

\begin{lema}\label{lagrange}
Let $\Cal L$ be a  loop with the inverse property and  suppose
that $\Cal H$ is a  finite index subloop of $\Cal L$ such that
$\Cal L$ has the Lagrange property relative to $\Cal H$. Let $\Cal
F$ be a subloop of $\Cal L$. Then $\Cal H\cap \Cal F$ is a finite
index subloop of $\Cal F$ and $\Cal F$ has  the  Lagrange property
relative to  $\Cal H\cap\Cal F$.
\end{lema}
\demo For all $a\in \Cal H\cap \Cal F$ and $f\in \Cal F,$ we have
$$(\Cal H\cap \Cal F)(af)=\Cal H(af)\cap \Cal F(af)=\Cal Hf\cap \Cal Ff=(H\cap \Cal
F)f,$$ thus  $\Cal F$ has the Lagrange property relative to $H\cap
\Cal F$. Let $n=[\Cal L:\Cal H]$ and $T=\{t_1,\dots t_n\}$ a
transversal of $\Cal H$. We can suppose that $\Cal Ht_j\cap \Cal
F\ne \emptyset$ if and only if $j\le m$ for some $m\le n$, thus
$\Cal F=\bigcup_{j=1}^m \Cal Ht_j\cap\Cal F.$ Let $s_j$ be an
arbitrary element of $\Cal Ht_j\cap\Cal F$. We claim that $(\Cal
H\cap\Cal F)s_j=\Cal Ht_j\cap\Cal F$. In fact, since $s_j=ht_j$
for some $h\in\Cal H$, then
$$(\Cal H\cap \Cal F)s_j=\Cal
H(ht_j)\cap \Cal Fs_j=\Cal Ht_j\cap \Cal
F,\quad\hbox{thus}\quad[\Cal H\cap \Cal F:\Cal F]=m.\fin$$

\begin{theorem} Let $R$ be a Dedekind algebraic integer domain, and
suppose that $R$ contains an infinite order unit. If $\Cal L$ is a
finite index subloop of $\GA$ and $\GA$ has the weak Lagrange
property, then $\Cal L$ is a congruence subloop. \end{theorem}
\demo Let $\Cal F$ be a finite index subloop of $\Cal H$, such
that $\GA$ has the Lagrange property relative to $\Cal F$. From
lemma \ref{lagrange} follows
$$[\GA_{(j)}:\GA_{(j)}\cap \Cal
F]\le [\GA:\Cal F]\le [\GA:\Cal L][\Cal L:\Cal F]<\infty.$$   Then
from Serre's Theorem follows that there are ideals $\q_j$ for
$j=1,2,3$, such that $\GA_{(j)}(\q_j)\subset \Cal F\cap
\GA_{(j)}\subset \Cal L$. Define $\q=\q_1\cap\q_2\cap\q_3$, then
$\GA_{(j)}(\q)\subset \Cal L$ for $j=1,2,3$, and from theorem
\ref{genera-n} follows that $\GA(\q)\subset \Cal F\subset \Cal
L$.\fin

\begin{proposicion} Let $R$ be a Dedekind algebraic integer domain
with an infinite order unit and $\q$ an ideal of $R$. Then
$\Delta(\q)=\GA(\q)$.
\end{proposicion}
\demo From Varser\v{s}ta\v{\i}n theorem $\Delta(\q)\cap \GA_{(j)}$
is a finite index subgroup of $\GA_{(j)}$, then from Serre's
Theorem there are $\q_j$ ideals such that $\GA_{(j)}(\q_j)\subset
\Delta(\q)\cap \GA_{(j)}\subset \Delta(\q)$. Thus, from theorem
\ref{genera-n} follows that $\GA(\q')\subset \Delta(\q)$ where
$\q'=\q_1\cap\q_2\cap\q_3$. Finally, from corollary \ref{Wolf}
follows that $\Delta(\q)=\GA(\q)$.\fin

\section{The  Loop  $SLL(2,\Z)$}

In this section, $\GA$ denotes  the loop $\SLL(2,\Z)$. For each
$n\ge 1$ $\GA(n)$ is a principal congruence subloop of level $n$
and
  $\GA'(n)$  is a subloop of  $\GA(n)$ generated by the associators and commutators,
i.e., the smallest loop that contains any element of the form
$$[A,B]=ABA^{-1}B^{-1}\qquad \hbox{and} \qquad
[A,B,C]=((AB)C)(A(BC))^{-1},$$ with $A,B$ and $C\in \GA(n)$. For
 $j=1,2,3,$ we denote $$T_j= \left[\begin{array}{cc} 0& e_j\\
-e_j
& 0\end{array}\right]\quad  S_j= \left[\begin{array}{cc} 1& e_j\\
0 & 1\end{array}\right]\quad\hbox{and}\quad
U_j=\left[\begin{array}{cc} 0& e_j\\ -e_j & 1\end{array}\right]$$
Observe that
 $\GA_{(j)}\cong SL(2,\Z)$   for each $j=1,2,3$. Specially,
$\GA_{(j)}$ is generated by two of the  following  matrices $T_j$,
$S_j$ and $U_j$ (See \cite{Ne} pag 139).

\begin{proposicion}\label{fin-ger} $\GA$ has minimal set of generators  $\{S_1,S_2,U_3\}$.
In general, for every positive integer $n$,  $\GA(n)$ is finitely
generated.
\end{proposicion}
\demo Since $T_1$, $T_2$, $S_3$ and $T_3$ can be written as a
product of $S_1$, $S_2$ and $U_3$ (see \cite{Voj} pag 190) and
$\GA$ is dissociative, it follows that this is a minimal set of
generator. For every $n>1$,  each $\GA_{(j)}(n)$ are finitely
generated free groups. The proposition follows as a trivial
consequences of theorem \ref{genera-n}.\fin

\begin{Observacion}\label{alfin} From proof of theorem \ref{genera-n}, we can
observe that $\GA(n)$ can be generated by the generators of
$\GA_{(1)}(n)$ and the matrices $\pmatrix{1&ne_j\cr 0 &1 }$ and
$\pmatrix{1&0\cr ne_j &1 }$ for $j=2,3$. Thus
$\GA(n)=\GA_{(1)}(n)\Delta(n)$.
\end{Observacion}



\begin{proposicion} Let  $p$ be a prime number  and $k\in\N.$
Then $\GA(p^k)$ is a normal subloop of  $\GA$ and it has index
$p^{7k}\left(1-\frac 1{p^4}\right)$.
\end{proposicion}
\demo Consider a loop homomorphism from $\GA$ onto
$\SLL(2,\Z_{p^k})$
$$\begin{array}{rcl}
\Theta:\GA&\longrightarrow& \SLL(2,\Z_{p^k})\\
A&\longmapsto& A\pmod {p^k}
\end{array}
$$
It is easy to prove that $\ker(\Theta)=\GA(p^k)$. Then
$$[\GA:\GA(p^k)]=\hbox{cardinality of } \SLL(2,\Z_{p^k}).$$
Now to obtain the cardinality of  $\SLL(2,\Z_{p^k})$, let us take
an arbitrary element $A=\pmatrix{a&(u_1,u_2,u_3) \cr
(v_1,v_2,v_3)&b}$, such that
$$\det A=ab-v_1u_1-v_2u_2-v_3u_3\equiv 1 \pmod {p^k}.$$
Observe  that $(a,v_1,v_2,v_3)$ can  assume any value different
from $p(n_1,n_2,n_3,n_4)$, i.e. this vector can  assume
$p^{4k}-p^{4(k-1)}$ different values. Fixing this vector, we know
that there is a coordinate non-divisible by  $p$. Without loss of
generality,  suppose that $a$ is not divisible by  $p$ (In the
case  $a$ is divisible by  $p$ there is some  $v_j$ that is not
divisible by  $p$ and  the argument follows  similarly). Then,
when we fix the values    $u_1$, $u_2$, $u_3$ the congruence
$$ab\equiv1+v_1u_1+v_2u_2+v_3u_3
 \pmod {p^k}$$ has an unique solution   $b$ modulo  $p^k$, i.e.,
 $u_1$, $u_2$ and $u_3$ determine exactly one value of  $b$
 modulo $p^k$ and thus $(b,u_1,u_2,u_3)$ can assume $p^{3k}$
 values. \fin

\begin{theorem}\label{indice} $\GA(n)$ is  a normal subloop of  $\GA$ with index
$\displaystyle{n^7\prod_{p|n}\left(1-\frac 1{p^4}\right)}.$
\end{theorem}
\demo Suppose  $n=p_1^{k_1}\cdots p_l^{k_l}$. Consider the
surjective loops homomorphism
$$\begin{array}{rcl}
\Theta:\GA&\longrightarrow& \prod\limits_{j=1}^l\SLL(2,\Z_{p_j^{k_j}})\\
A&\longmapsto& \prod\limits_{j=1}^l A\pmod {p_j^{k_j}}
\end{array}
$$
It is easy to prove that $\ker(\Theta)=\GA(n)$ and then  $\GA(n)$
is a normal subloop of $\GA$, in addition
$$[\GA:\GA(n)]=\hbox{number of elements of } \prod_{j=1}^l\SLL(2,\Z_{p_j^{k_j}}).$$
Then the  theorem follows from the   proposition before.\fin

\begin{lema}\label{mcd} Let  $\Cal L$ be a normal subloop  $\GA$, $m$ and $n$
positive integers such that  $\Delta(n)\subset \Cal L$ and
$\Delta(m)\subset \Cal L$. If $d=(n,m)$, then $\Delta(d)\subset
\Cal L$.
\end{lema}
\demo Since  $d=tn+sm$ where  $t,s\in\Z$, then lemma follows from
$$\left[\begin
{array}{cc}1 & dx \\
0 & 1\end{array}\right]=\left[\begin
{array}{cc}1 & nx \\
0 & 1\end{array}\right]^t\left[\begin
{array}{cc}1 & mx \\
0 & 1\end{array}\right]^s\quad\hbox{and}\quad\left[\begin
{array}{cc}1 & 0 \\
dx & 1\end{array}\right]=\left[\begin
{array}{cc}1 & 0 \\
nx & 1\end{array}\right]^t\left[\begin
{array}{cc}1 & 0 \\
mx & 1\end{array}\right]^s$$ for every $x\in \Z^3$.
 \fin

 It is known that
 $\Cal \GA'(n)$
 is a normal subloop of $\GA(n)$ (ver \cite{GJM} p\'ag. 56).
Let $\G(n)$ denote the group $\frac {\GA(n)}{\GA'(n)}$. From
proposition \ref{fin-ger} follows that $\GA(n)$ is finitely
generated, and thus $\G(n)$ is a finite generated abelian group.

\begin{lema} For $n>5$,  $\G(n)$ is infinite.
\end{lema}
\demo For each $A\in \GA(n)$, from observation \ref{alfin} we know
that $A=BC$ where $B\in \GA_{(1)}(n)$ and $C\in \Delta(n)$. Let
$\Theta:\GA(n)\to \frac {\GA_{(1)}(n)}{\Delta_{(1)}(n)}$ defined
by the rule $\Theta(A)=B\Delta_{(1)}(n)$. To show that $\Theta$ is
well defined, suppose that $A=B_1C_1=B_2C_2$ where $B_1,B_2\in
\GA_{(1)}(n)$ and $C_1,C_2\in \Delta(n)$. Since $\Delta(n)$ is
normal,
$$B_1\Delta(n)=(B_1C_1)\Delta(n)=(B_2C_2)\Delta(n)=B_2\Delta(n),$$
thus $B_1^{-1}B_2\in \Delta(n)\cap \GA_{(1)}(n)=\Delta_{(1)}(n)$
and  $\Theta$ is a loop homomorphism with kernel $\Delta(n)$. It
follows that
$$\frac {\GA(n)}{\Delta(n)}\simeq \frac
{\GA_{(1)}(n)}{\Delta_{(1)}(n)}.$$ Now, from the second and third
homomorphism theorems  for loops, we have
$$\frac {\frac {\GA(n)}{\Delta(n)}}{\left(\frac
{\GA(n)}{\Delta(n)}\right)'}\simeq \frac {\frac
{\GA(n)}{\Delta(n)}}{\frac{\GA'(n)}{\GA'(n)\cap\Delta(n)}}\simeq
\frac {\frac
{\GA(n)}{\Delta(n)}}{\frac{\GA'(n)\Delta(n)}{\Delta(n)}}\simeq
\frac {\GA(n)}{\GA'(n)\Delta(n)}.$$ Let $C(n)$ denote the group
$\frac {\GA_{(1)}(n)}{\Delta_{(1)}(n)}$. In  \cite{BLS} it  proved
that profinite cohomology group
$$H^1(\lim\limits_{\leftarrow}
C(n),\Q/\Z)=Hom(\frac{C(n)}{(C(n))'}, \Q/\Z)=0$$ and
$$H^1(\lim\limits_{\leftarrow}
C(n),\Q/\Z)=Hom(\frac{C(n)}{(C(n))'}, \Z)$$ is infinity,  thus
$\frac{C(n)}{(C(n))'}$ is a abelian torsion-free
 group  for
all $n\ge 1$. Since $C(n)$ is an   infinite group for $n>5$ (see
\cite{Ne} pag 145),  it follows that $[C(n):(C(n))']$ is infinite
for $n>5$ and
$$|\G(n)|=[\GA(n):\GA'(n)]\ge
[\GA(n):\GA'(n)\Delta(n)]=[C(n):(C(n))']=\infty.\fin$$

Denote $\G_s(n)$ the subgroup of $\G(n)$ generated by the $s$-th
powers, i.e. $\G_s(n)=\langle A^s| A\in \G(n)\rangle$, then
$\G(n)/\G_s(n)$ is finite, in fact,  $[{\G(n)}:{\G_s(n)}]\le s^k$
where $k$ is the number of generators of  $\G(n)$.  Observe that
the homomorphisms
$$\GA(n)\stackrel{\pi}{\longrightarrow}
\G(n)\stackrel{\psi}{\longrightarrow}\frac {\G(n)}{\G_s(n)}$$ are
well defined and surjective. Denote $\GA(n,s)=\ker (\psi\circ
\pi)=\pi^{-1}(\G_s(n))$. It is easy to see that $\GA(n,s)$ is
generated by the commutators, associators  and  the set
$\{A^s|A\in \GA(n)\}$.

\begin{theorem} Let  $n>5$   and $s$ be an odd integer   such that $(n,s)=1$. Then $\GA(n,s)$ is
a finite index subloop of  $\GA$ that is not a congruence subloop.
\end{theorem}
\demo Since $$[\GA:\GA(n,s)]=[\GA:\GA(n)][\GA(n):\GA(n,s)],$$ from
 theorem  \ref{indice}, the definition of  $\GA(n,s)$ and the
homomorphisms theorem we have
$$[\GA:\GA(n)]=n^7\prod_{p|n} \Bigl(1-\frac 1{p^4}\Bigr)\quad\hbox{and}\quad
[\GA(n):\GA(n,s)]=[{\G(n)}:{\G_s(n)}]$$ thus  $[\GA:\GA(n,s)]$ is
finite. In addition, since $\G(n)$ is a finite generated abelian
group and for $n>5$  is infinite, it follows that $\G(n)$ has a
factor isomorphic to $\Z$, and thus
$$[\GA(n):\GA(n,s)]=[{\G(n)}:{\G_s(n)}]\ge s>1,$$ in particular
$\GA(n,s)\varsubsetneq \GA(n)$. Now, suppose $\GA(n,s)$ is a
congruence subloop, since $A_{jn}^{2n},B_{jn}^{2n}\in \GA'(n)$ for
every  $j=1,2,3$, then
$$\Delta(2n^2)\subset \GA'(n)\subset \GA(n,s).$$ In the same way,
from the definition of  $\GA(n,s)$, we have that $A^s\in \GA(n,s)$
for every $A\in \GA(n)$, in particular $\pmatrix{ 1&nx\cr 0&1}^s,
\pmatrix{ 1&0\cr nx&1}^s \in \GA(n,s)$ for every $x\in \Z^3$. Then
$\Delta(ns)\subset \GA(n,s)$. Now, $(ns,2n^2)= n$, then   from
lemma \ref{mcd}, we have  $\Delta(n)\subset \GA(n,s)$ and, finally
from corollary \ref{Wolf}, it follows that
$$\GA(n)\subset \GA(n,s)\varsubsetneq \GA(n),$$  but this is impossible. \fin

\noindent If $s>120$ is even, the theorem above is also true. In
fact, since $(ns,2n^2)=2n$, from lemma \ref{mcd} we have
$\Delta(2n)\subset \GA(n,s)$,  it follows that  $\GA(2n)\subset
\GA(n,s)$ and
$$s\le [\GA(n):\GA(n,s)]\le [\GA(n):\GA(2n)]= 120,$$
but this is a contradiction, since $s>120$.

\begin{center}
\textsc{Acknowledgmet}
\end{center}

The authors are grateful to Guilherme Leal for suggesting the
problem and helpful conversations.

\end{document}